# Fuzzy location and allocation Hub Network Design for Air Cargo Transportation Considering Sustainability and Time Window

Ali Mohammad Malekdar[a], Mohsen Akbarpour Shirazi[b]

[a,b] Department of Industrial Engineering & Management Systems, Amirkabir University of Technology (Tehran Polytechnic)

**Corresponding Author**: Mohsen Akbarpour Shirazi
**Email**: Akbarpour@aut.ac.ir

**Abstract:**
Hub location Problems seek to find hub facilities and assign non-hub nodes to them in such a way that the flow between origin and destination should be effectively established according to the desired goal. In general, in the literature of location, it is assumed that the time horizon of hub network design is a single time horizon. In the last two decades these problems have attracted special attention in the field of facility location problems and have wide applications in different fields including air cargo transportation. Cargo transportation is one of the most important economic sectors of any country. There are different ways to transport cargo, but air transport is preferred because it has high speed and security, so it is suitable for transporting goods related to technology, food, medicines, etc. In this article designing a hub network for air cargo transportation, taking into account hard and soft time windows along with considering the limited capacity for each hub under uncertainty is discussed. The proposed model is a developed model of an existing model in literature. Our Study has three linear functions: economic, environmental and social. In this article method of Fuzzy programming has been used to control the non-deterministic demand parameter. Results of model that has been solved by epsilon limitation method, NSGA-II, MOPSO and MOWOA algorithm show that as the uncertainty rate increases, the total costs of the system as well as the amount of environmental pollution increases. The reviews indicate the high performance of the NSGA-II algorithm in solving the proposed model.

**Introduction**

Air transportation is one of the most complex systems created by humans. It is one of the most important methods for the rapid transfer of valuable goods. However, air transportation has always had a much smaller share compared to other modes of cargo transport. Essentially, air transportation has always been more favored by passengers, and the majority of the industry's revenue comes from passenger transportation. According to statistics, 75% of airline revenue comes from passenger transport, 15% from cargo transport, and 10% from other sources (Shen et al., 2021). In contrast, other transportation industries such as ground and maritime transport primarily earn their revenue from cargo transportation. However, it is expected that due to the increasing global trade, the need for transporting goods over longer distances, and the booming industry of lightweight but valuable electronic components like microchips, the air transportation industry is growing. The fundamental reason for the limited use of air transportation for cargo is its high cost. The strength of air transportation compared to other methods is its high speed for long-distance transfers. Therefore, if time is a crucial factor in the transport of cargo, producers are willing to bear the high cost (Alumur et al., 2021). For example, technology owners often use air transportation because the high cost is offset by the benefit of getting their products to market earlier. Air transportation is also used for perishable goods such as fresh food, pharmaceuticals, and agricultural products.

Airlines use two methods for transporting cargo to destinations: dedicated cargo aircraft and passenger aircraft. The second method occurs when an airline decides to use the extra cargo space in passenger aircraft to transport smaller cargo shipments (Delgado et al., 2020). Each of these methods has its own advantages and dis advantages. For example, using extra space in passenger aircraft is usually less costly for shippers and airlines. On the other hand, dedicated cargo aircraft offer more flexibility due to their specialized equipment and can handle larger and more sensitive shipments. Additionally, the origin, destination, and schedule of "all-cargo" flights are precisely planned to match the demands of cargo brokers and other shippers, which is another significant advantage of using cargo aircraft over passenger aircraft. As a result, dedicated cargo aircraft typically provide higher-quality service and greater safety. The importance of air transport has led researchers to use hub location problems for this purpose. In hub location problems, nodes are designated as hubs to reduce direct cargo transport costs, and other nodes are allocated to these hubs. Essentially, cargo distribution occurs between the hubs. Recent research has increasingly focused on sustainability in hub location within air transport systems, designing sustainable hub networks considering environmental, economic, and social aspects. The goal is to balance the various objectives of a supply chain. Traditional hub network design models aim to minimize costs or, in other words, maximize profit. Their objective is to minimize the total cost of the hub network, including fixed hub establishment costs, transportation costs, and other expenses. Some researchers have developed bi-objective models with economic-environmental objective functions. Their aim is to reduce environmental impacts (Pourghader et al., 2021).

According to relevant reports, greenhouse gas emissions from air transport are increasing at a rate of 3.6% per year. The fuel burned by aircraft engines releases pollutants such as carbon dioxide ($CO_2$), carbon monoxide (CO), water vapor ($H_2O$), hydrocarbons (HC), nitrogen oxides (NOX), and sulfur dioxide, which are key environmental impacts on the atmosphere (Krile et al., 2015). The concept of "green aviation" represents new conditions in the air transport industry for environmental protection, playing a crucial role in the sustainable development of airlines. Therefore, comprehensive attention to sustainability concepts in all aspects is essential for hub networks, which can contribute to sustainable development. It is estimated that by 2040, CO2 emissions from air transport will increase by 21%, and NOX emissions by 16% (Harley et al., 2020). Research indicates that among various modes of transport, such as air, rail, road, and water, air transport has the highest share in the increase of greenhouse gas emissions. Reducing emissions from the aviation industry might be an effective climate policy (Zhang et al., 2016).

Given that many hub location studies focus on a single objective and cost reduction, service level objectives like delivery time are often overlooked. Cost and delivery time are conflicting objectives often ignored in many service networks. Hub location problems such as P-Hub center location problem and P-Hub covering location problem focus specifically on service level objectives. P-Hub center center problems aim to minimize the maximum service time between origin-destination pairs (Ernst et al., 2009). On the other hand, in P-Hub covering problems, demand between an origin-destination pair is covered if it can be met within a certain time (Campbell, 1994). Meanwhile, cargo transport companies are eager to minimize both fixed and operational costs of establishing hub networks and total transportation costs to provide services within promised service times (Boysen and Fliedner, 2010). Considering real-world conditions, reducing transportation costs is only one part of the decision-making criteria. Service quality must also be considered in hub network design decisions. For instance, in the cargo shipment sector, a service time window is usually set for origin-destination pairs. If service time is ignored in hub network design, the resulting hub network may not be feasible in terms of service level in the future.

Passenger airline networks have been a major application area for hub location models for years. Since the demand in these systems consists of individuals who do not want to wait, most real-world airline networks use a multiple allocation network. Each origin-destination pair is then served by the fastest or cheapest route through the hub. Due to the greater time flexibility for cargo transport in the aviation industry, it allows for reduced costs and ease in managing transfers. Therefore, as long as service levels are maintained, there is flexibility in route selection. This capability allows transport companies to use routes with more hubs, which might increase distance or travel time (compared to direct routes) to achieve cost savings. Air transport, due to long distances and the amount of cargo carried, usually has a significant share in greenhouse gas emissions. In such conditions, selecting cities as hubs and determining network routes in a way that minimizes transportation costs within the network and the amount of targeted gas emissions is highly beneficial. This article presents a model for the capacitated hub location and allocation problem in the air freight transportation industry, aiming to minimize network transportation costs, reduce emitted pollutants, and minimize penalties assigned for failing to meet the specified time window.

The structure of the paper is as follows. Section two reviews the research literature and examines studies on hub location in the aviation industry. Section three first designs an uncertain model for hub location in the air transportation transport system and then presents a fuzzy programming method to control demand parameters. Section four introduces solution methods such as epsilon constraint method and MOPSO, MOWOA, and NSGA-II algorithms. Section five analyzes various numerical examples and sensitivity analysis and prioritizes solution methods using TOPSIS. Section six concludes with findings and suggestions for future research.

**Literature Review**

Research in the field of hub location has been extensively studied by many researchers. This section reviews the most important articles on hub location in air or ground transportation systems. Yaman et al. (2007) modeled a hub location problem for cargo delivery in Turkey. The computational results based on the developed model indicated improvements in Turkey's cargo transportation fleet. Lin et al. (2012) presented a median hub location model to minimize distribution costs of air transportation in China. They used a genetic algorithm to find the best solution. Sensitivity analysis on discount factor showed that economies of scale in the main routes of hub-and-spoke networks could significantly impact operational costs and routing patterns. Ambrosino and Sciomachen (2016) investigated a capacitated hub location problem for a multi-level cargo transportation network. Their goal was to minimize hub location and transportation costs in the network. They used CPLEX to solve problems of various sizes.

Shang et al. (2020) proposed a multi-modes hierarchical stochastic hub location model for cargo delivery systems. The model presented in this paper is of the star model. They employed a genetic algorithm to solve the problem. Computational results showed that with increasing levels of confidence, airport hubs in the cargo distribution network are located further apart to gain more time advantages. Golestani et al. (2021) proposed a bi-objective green hub location problem aiming to minimize total costs (including transportation, hub establishment, storage temperature regulation, and carbon emissions) and maximize the quality of delivered products to the customers. They utilized the epsilon-constraint method to solve the bi-objective model. They demonstrated that distributing perishable products with different storage temperatures in the proposed model could help maintain the quality of delivered items while reducing overall system costs and considering carbon emissions.

Shang et al. (2021) provided a heuristic algorithm for the bi-objective multi-modes hierarchical hub location problem for cargo delivery systems. Their objectives were to minimize total system costs and the maximum delivery time. They used the epsilon-constraint method and NSGA-II, showing that NSGA-II performs better than other methods. Zhu et al. (2023) examined a hub location problem for air transportation distribution during the COVID-19 pandemic. In this mathematical model, the demand parameter was considered as an uncertain parameter. Shipments are placed in air containers based on weight and volume, then flown from regional collection points to a hub for integration before transport to subsequent destinations. They used a genetic algorithm to minimize total costs. Eydi and Shirinbayan (2023) modeled a hierarchical multi-commodity hub location problem with fuzzy demand. Their main goal was to minimize total transportation costs in the network to determine optimal hub locations, assign non-hub nodes to hubs, and identify the type of vehicles needed for each. They proposed using a genetic algorithm to solve the problem.

Roozkhosh et al. (2023) presented a new model for the hub location-allocation problem, which can evaluate the impact of costs and delays in delivery to destinations, find the optimal number of hubs and vehicles with different capacities, and address congestion in hubs. They employed a particle swarm optimization algorithm (PSO) to solve the problem using data from Australia Post. Rahmati et al. (2024)

examined a hub location problem with the goal of maximizing random profit in a two-stage stochastic model with uncertain demand. They used the advanced sample average approximation (ASAA) method to determine the appropriate number of scenarios. Computational experiments showed that all carbon regulations could reduce overall carbon emissions. Andaryan et al. (2024) studied a single hub location problem with Bernoulli demands and used the Tabu search algorithm as a solution approach to handle large problem instances.

The literature review indicates that most articles have addressed only the economic and environmental aspects of the problem, with the social aspects of air transportation distribution systems not being comprehensively studied. Therefore, this paper addresses this research gap by presenting a comprehensive model of sustainable hub location.

While previous studies have focused on cost minimization or environmental impacts in hub network design, few have addressed the comprehensive integration of sustainability, time window constraints, and uncertainty. For example, studies by Golestani et al. (2021) and Shang et al. (2021) primarily focus on bi-objective models that do not consider the social dimension of sustainability or the uncertainty in demand. In contrast, this study presents a three-objective model that simultaneously optimizes economic, environmental, and social outcomes, while also incorporating fuzzy demand to reflect real-world uncertainties.

Moreover, the proposed use of metaheuristic algorithms such as NSGA-II and MOWOA represents a significant advancement over the traditional exact methods used in many previous works, enabling the model to handle larger and more complex problems efficiently.

**Statement of the Problem and Modeling**

Despite the significant importance of air cargo transportation in modern logistics systems, many existing studies on hub network design primarily focus on minimizing costs, often neglecting the critical aspects of sustainability and time-sensitive logistics. Furthermore, most models assume a deterministic environment, which does not accurately reflect real-world uncertainties such as fluctuating demand. This gap in the literature, especially in the context of incorporating both fuzzy demand and time window constraints, highlights the need for more robust models that address not only economic factors but also environmental and social objectives. This study addresses these gaps by proposing a comprehensive hub network design model for air cargo transportation that integrates sustainability goals, uncertainty in demand, and hard and soft time windows. This novel approach ensures that the model is more applicable to real-world logistics, where uncertainties and sustainability are becoming increasingly crucial.

This study addresses these gaps by proposing a comprehensive hub network design model for air cargo transportation that integrates sustainability goals, uncertainty in demand, and hard and soft time windows. This novel approach ensures that the model is more applicable to real-world logistics, where uncertainties and sustainability are becoming increasingly crucial.

The literature review highlights the importance of air freight transportation and the use of hub location in distribution network problems. This section presents a model for the hub location problem for air freight transportation. In this model, there is a set of nodes where the distribution of goods between any two nodes is carried out by the air transportation system. In this problem, the number of nodes that can be converted into hubs is specified, and the goal is the optimal allocation of non-hub nodes to hubs and the distribution of goods via the air transportation system. Considering that the aim of presenting the mathematical model is not merely to reduce costs, environmental and social aspects are also taken into account in this mathematical model.

In the past decade, hub network design has attracted much attention in operational research. Moreover, sustainability plays a crucial role in the performance of hub networks. In this study, a method used by the European Environment Agency (EEA, 2016) for estimating aircraft emissions is employed. In this method, the mechanism and operational cycle of aircraft, as shown in Fig. 1, are divided into two main parts:

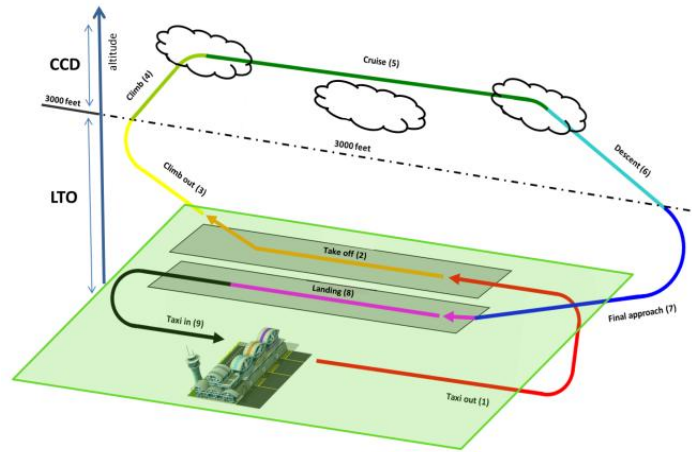

**Figure 1- Operating cycle of the aircraft during the flight (Intergovernmental Panel on Climate Change (IPCC), 2000)**

- LTO Cycle (The Landing/Take-Off): This cycle includes all activities that occur near airports at altitudes below 3,000 feet during the departure and arrival phases of the flight.
- CCD Cycle (Climb/Cruise/Descent): The cruise cycle includes all activities that occur at altitudes above 3,000 feet during flight.

In 2016, the EEA published several methods called "Tier 3" for calculating aircraft emissions. According to this method, the amount of pollutants ($VOC, NO_2, CO$) emitted over a $d_{ij}$ distance is calculated as $E_{p_1}^{LTO} + R_{p_1}^{CCD}(d_{ij})$. Here, $E_{p_1}^{LTO}$ represents the amount of pollutants ($VOC, NO_2, CO$) emitted during the LTO cycle, and $R_{p_1}^{CCD}(d_{ij})$ represents the amount of pollutants ($VOC, NO_2, CO$) emitted during the cruise cycle. Additionally, the amount of pollutants ($SO_2, CO_2$) emitted over a $d_{ij}$ distance is calculated as $E_{p_2}^{LTO} + R_{p_2}^{CCD}(d_{ij})$. Here, $E_{p_2}^{LTO}$ represents the amount of pollutants ($SO_2, CO_2$) emitted during the LTO cycle, and $R_{p_2}^{CCD}(d_{ij})$ represents the amount of pollutants ($SO_2, CO_2$) emitted during the cruise cycle.

On the other hand, in real-world transportation problems, customers prefer to receive timely scheduling for the delivery of goods or services. However, it is important to note that in practice, there may be situations where some service providers face problems and violate the schedule to which they were committed. It should be noted that if a constraint must be satisfied, it is called "hard," whereas if it can be violated, it is called "soft." Violating soft constraints is usually penalized and then added to the objective function. In this mechanism, the time between the first and last requested service times is referred to as the requested time window. Demands that are fulfilled within this time window incur no penalties. Otherwise, they face penalties in the objective function. The goal in this case of sustainability is to design the network and allocate in such a way that it focuses on timely delivery.

As shown in Fig. 2, in the proposed soft time window, $E_{ij}$ represents the lower bound of the time window, and $L_{ij}$ represents the upper bound of the time window. Ideally, if the demand is met within the specified time, which falls between $E_{ij}$ and $L_{ij}$, there will be no penalty. Otherwise, a specific penalty will be incurred.

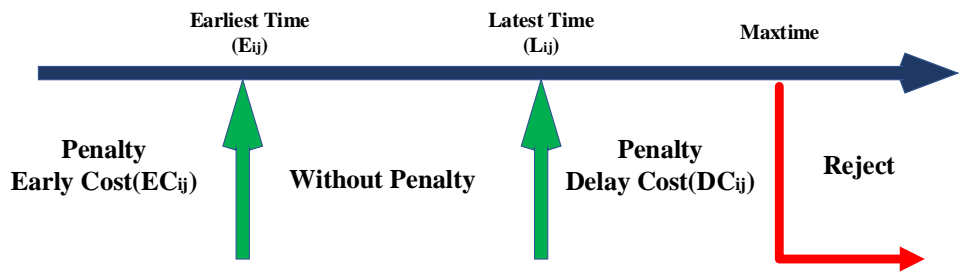

**Figure 2- Time window for delivery of goods in the air transport system**

Thus, the problem has three objectives: the first objective minimizes transportation costs within the network, the second objective minimizes the emission of pollutants, and the third objective ensures that shipments arrive within the requested time window.

Based on these objectives, the assumptions of the mathematical model are as follows:

- There is a set of nodes that serve as origins and destinations of flows.
- A non-hub node is connected to a hub node.
- A hub node can be connected to multiple non-hub nodes.
- The speed and capacity of the air fleet are uniform.
- All nodes have sufficient flights to meet demands.
- Hub nodes have limited capacity.
- There is a time window for sending shipments from origin to destination.
- Demand is considered an uncertain parameter and modeled as a trapezoidal fuzzy number.
- The number of hub locations is predetermined and fixed.
- Various types of connections between origin and destination nodes are possible.

Given the available objectives in the mathematical model and the stated assumptions, the symbols used in the mathematical model are defined as follows:

**Parameters:**

$N$  Total number of nodes $i, j, k, l \in N$

$p$  Maximum number of predetermined hubs

$\omega$  Maximum allowable distance between a non-hub node and the assigned hub node

$f_k$  Fixed operational cost of a hub $k \in N$

$\vartheta_k$  Capacity of goods handling at the hub $k \in N$

$u_k$  Cost of handling each unit of cargo at the hub $k \in N$

$\tau_{ij}$  Delivery time of goods between nodes $i \in N$ and $j \in N$

$d_{ij}$  Distance between nodes $i \in N$ and $j \in N$

$\sigma_{ij}$  Maximum transfer time between two nodes $i \in N$ and $j \in N$

$c_{ij}$  Transportation cost per unit of cargo between nodes $i \in N$ and $j \in N$

$\tilde{q}_{ij}$  Uncertain demand for cargo between nodes

$\alpha$  Discount factor for transportation costs (from hub to hub)

$\beta$    Discount factor for transportation costs (from non-hub to hub and vice versa)

$EC_{ij}$    Penalty cost for meeting demand before the time window between nodes $i \in N$ and $j \in N$

$DC_{ij}$    Penalty cost for meeting demand after the time window between nodes $i \in N$ and $j \in N$

$E_{ij}$    Lower bound of the time window for sending shipments between nodes $i \in N$ and $j \in N$

$L_{ij}$    Upper bound of the time window for sending shipments between nodes $i \in N$ and $j \in N$

$\varphi$    Capacity of the aircraft

$E_{p_1}^{LTO}$    The amount of emission of pollutants ($VOC, NO_2, CO$) emitted during the LTO cycle

$E_{p_2}^{LTO}$    The amount of emission of pollutants $SO_2, CO_2$ during the LTO cycle

$R_{p_1}^{CCD}(d_{ij})$    Emission of pollutants $VOC, NO_2, CO$ during the cruise cycle

$R_{p_2}^{CCD}(d_{ij})$    Emission of pollutants $SO_2, CO_2$ during the cruise cycle

**Decision Variables:**

$X_k$ 1 if node $k \in N$ is selected as a hub, 0 otherwise.

$Z_{ik}$ 1 if non-hub node $i \in N$ is assigned to hub node $k \in N$, 0 otherwise.

$Y_{ij}$ 1 if direct transport is established between origin node $i \in N$ and destination node $j \in N$, 0 otherwise.

$W_{ikj}$ 1 if transport between origin node $i \in N$ and destination node $j \in N$ is through a $k \in N$ hub, 0 otherwise.

$V_{iklj}$ 1 if transport between origin node $i \in N$ and destination node $j \in N$ is through $k \in N$ hub and then through hub $l \in N$, 0 otherwise.

Based on the defined symbols, the multi-objective mathematical model for hub location in air transportation under uncertainty is as follows:

$$\begin{aligned} Min\ Z_1 &= \sum_{i \in N} \sum_{j \in N} c_{ij} d_{ij} \tilde{q}_{ij} Y_{ij} \\ &+ \sum_{i \in N} \sum_{j \in N} \sum_{k \in N} \beta(c_{ik} d_{ik} + c_{kj} d_{kj}) \tilde{q}_{ij} W_{ikj} \\ &+ \sum_{i \in N} \sum_{j \in N} \sum_{k \in N} \sum_{l \in N} u_l \tilde{q}_{ij} V_{iklj} + \\ &\sum_{i \in N} \sum_{j \in N} \sum_{k \in N} \sum_{l \in N} \left( \beta(c_{ik} d_{ik} + c_{kj} d_{kj}) + \alpha(c_{kl} d_{kl}) \right) \tilde{q}_{ij} V_{iklj} \\ &+ \sum_{i \in N} \sum_{j \in N} \sum_{k \in N} u_k \tilde{q}_{ij} W_{ikj} \\ &+ \sum_{k \in N} f_k X_k \end{aligned} \quad (1)$$

$$\begin{aligned} Min\ Z_2 &= \sum_{i \in N} \sum_{j \in N} \left( E_{p_1}^{LTO} + R_{p_1}^{CCD}(d_{ij}) \right) \left\lceil \frac{\tilde{q}_{ij}}{\varphi} \right\rceil Y_{ij} \\ &+ \sum_{i \in N} \sum_{j \in N} \sum_{k \in N} \left( 2 E_{p_1}^{LTO} + R_{p_1}^{CCD}(d_{ik} + d_{kj}) \right) \left\lceil \frac{\tilde{q}_{ij}}{\varphi} \right\rceil W_{ikj} + \\ &\sum_{i \in N} \sum_{j \in N} \sum_{k \in N} \sum_{l \in N} \left( 3 E_{p_1}^{LTO} + R_{p_1}^{CCD}(d_{ik} + d_{kl} + d_{kj}) \right) \left\lceil \frac{\tilde{q}_{ij}}{\varphi} \right\rceil V_{iklj} \\ &+ \sum_{i \in N} \sum_{j \in N} \left( E_{p_2}^{LTO} + R_{p_2}^{CCD}(d_{ij}) \right) \left\lceil \frac{\tilde{q}_{ij}}{\varphi} \right\rceil Y_{ij} + \\ &\sum_{i \in N} \sum_{j \in N} \sum_{k \in N} \left( 2 E_{p_2}^{LTO} + R_{p_2}^{CCD}(d_{ik} + d_{kj}) \right) \left\lceil \frac{\tilde{q}_{ij}}{\varphi} \right\rceil W_{ikj} + \\ &\sum_{i \in N} \sum_{j \in N} \sum_{k \in N} \sum_{l \in N} \left( 3 E_{p_2}^{LTO} + R_{p_2}^{CCD}(d_{ik} + d_{kl} + d_{kj}) \right) \left\lceil \frac{\tilde{q}_{ij}}{\varphi} \right\rceil V_{iklj} \end{aligned} \quad (2)$$

$$\begin{aligned} Min\ Z_3 &= \sum_{i \in N} \sum_{j \in N} EC_{ij} \max\{0, E_{ij} - \tau_{ij}\} Y_{ij} \\ &+ \sum_{i \in N} \sum_{j \in N} \sum_{k \in N} EC_{ij} \max\{0, E_{ij} - \tau_{ik} - \tau_{kj}\} W_{ikj} + \end{aligned} \quad (3)$$

$$\sum_{i \in N}\sum_{j \in N}\sum_{k \in N}\sum_{l \in N} EC_{ij} \max\{0, E_{ij} - \tau_{ik} - \tau_{kl} - \tau_{lj}\} V_{iklj}$$

$$+ \sum_{i \in N}\sum_{j \in N} DC_{ij} \max\{0, \tau_{ij} - L_{ij}\} Y_{ij} +$$

$$\sum_{i \in N}\sum_{j \in N}\sum_{k \in N} DC_{ij} \max\{0, \tau_{ik} + \tau_{kj} - L_{ij}\} W_{ikj}$$

$$+ \sum_{i \in N}\sum_{j \in N}\sum_{k \in N}\sum_{l \in N} DC_{ij} \max\{0, \tau_{ik} + \tau_{kl} + \tau_{lj} - L_{ij}\} V_{iklj}$$

s.t.:

$$Z_{ik} \leq X_k, \quad \forall i, k \in N \quad (4)$$

$$\sum_{k \in N} Z_{ik} = 1, \quad \forall i \in N \quad (5)$$

$$\sum_{k \in N} X_k \leq p \quad (6)$$

$$d_{ij} Z_{ij} \leq \omega, \quad \forall i, j \in N \quad (7)$$

$$W_{ikj} \leq Z_{ik}, \quad \forall i, j, k \in N \quad (8)$$

$$W_{ikj} \leq Z_{jk}, \quad \forall i, j, k \in N \quad (9)$$

$$V_{iklj} \leq Z_{ik}, \quad \forall i, j, k, l \in N \quad (10)$$

$$V_{iklj} \leq Z_{kl}, \quad \forall i, j, k, l \in N \quad (11)$$

$$V_{iklj} \leq Z_{lj}, \quad \forall i, j, k, l \in N \quad (12)$$

$$Y_{ij} + \sum_{k \in N} W_{ikj} + \sum_{k \in N}\sum_{l \in N} V_{iklj} = 1, \quad \forall i, j \in N \quad (13)$$

$$\sum_{i \in N}\sum_{j \in N} \tilde{q}_{ij} W_{ikj} + \sum_{i \in N}\sum_{j \in N}\sum_{l \in N} \tilde{q}_{ij} V_{iklj} + \sum_{i \in N}\sum_{j \in N}\sum_{l \in N} \tilde{q}_{ij} V_{ilkj} \leq \vartheta_k + M(1 - X_k), \quad \forall k \in N \quad (14)$$

$$\tau_{ij} Y_{ij} + (\tau_{ik} + \tau_{kj}) W_{ikj} + (\tau_{ik} + \tau_{kl} + \tau_{lj}) V_{ilkj} \leq \sigma_{ij}, \quad \forall i, j, k, l \in N \quad (15)$$

$$X_k, Z_{ik}, Y_{ij}, W_{ikj}, V_{ilkj} \in \{0,1\}, \quad \forall i, j, k, l \in N \quad (16)$$

Equation 1 shows the objective function of minimizing the total costs. The objective function of the total cost includes six components: 1. Transportation costs of air transportation (direct connection), 2. Transportation costs of air transportation (connection via one hub), 3. Transportation costs of air transportation (connection via two hubs), 4. Handling costs within the hub (connection via one hub), 5. Handling costs within the hub (connection via two hubs), 6. Operational costs of establishing hubs.

Equation (2) represents the objective function for minimizing aircraft emissions. In this equation, the mechanism and operational cycle of the aircraft are divided into two general parts: the LTO cycle and the CCD cycle.

Equation (3) shows the social objective function of the problem, which includes timely delivery of goods to customers. In this equation, the aim is to minimize the penalty resulting from the time discrepancy of shipments for each origin-destination pair.

Equation (4) ensures that a non-hub node can only be connected to a hub node where a hub has been established.

Equation (5) states that a non-hub node can be connected to at most one hub node.

Equation (6) guarantees that the number of hubs established in the network does not exceed a predetermined value.

Equation (7) indicates that a non-hub node can only be connected to a hub node if the distance between them is less than a predetermined permissible limit.

Equation (8) and (9) state that when two nodes are connected by a hub, both non-hub nodes should be assigned to the same hub node.

Equation (10) to (12) indicate that in the connection of two nodes via two hub nodes, the non-hub node $i \in N$ should be connected to the hub node $k \in N$, and the non-hub node $j \in N$ should be connected to the hub node $l \in N$, and the hub node $k \in N$ should be connected to the hub node $l \in N$.

Equation (13) states that only one connection is allowed between each origin-destination pair (direct connection, connection via one hub, or connection via two hubs).

Equation (14) indicates that each hub node has a limited capacity.

Equation (15) ensures that the shipping time for each origin-destination pair does not exceed a predetermined limit.

Equation (16) defines the type of decision variables.

As stated in the assumptions of the mathematical model, the demand parameter is considered as an uncertain parameter and is represented by a trapezoidal fuzzy number. Therefore, to control this parameter, the fuzzy programming method has been used. Consider the following linear mathematical programming model with the fuzzy demand parameter:

$$Min\ Z = c^t x \quad (17)$$

$s.t.:$

$$x \in N(\tilde{A}, \tilde{B}) = \{x \in R^n | a_{ij}x \geq \widetilde{q_{ij}}\}, \quad i \in m, j \quad x \geq 0 \quad (18)$$

where the parameter $c = (c_1, c_2, \ldots, c_n), A = [a_{ij}]_{m \times n}, \widetilde{q_{ij}} = (\tilde{q}_{11}, \tilde{q}_{12}, \ldots, \tilde{q}_{mn})^t$ used in the objective function is the coefficient vector, and the right-hand side is the parameter of the constraint (the demand value). The probability distribution function of the fuzzy demand parameter is assumed based on the characteristics of fuzzy numbers.

## Solution Methods

Given the multi-objective nature of the mathematical model presented in this paper, the epsilon-constrained method is used for validation and sensitivity analysis, while algorithms such as NSGA-II, MOPSO, and MOWOA are employed to solve larger instances of the model. This section describes these solution methods in detail.

### Epsilon-Constrained Method

In the epsilon-constrained method, one objective function is selected for optimization while the other objective functions are considered as constraints. These constraints have upper bounds defined by small values called epsilon. This

Finally, $x = (x_1, x_2, \ldots, x_n)$ represents the decision vector. For the feasibility and optimization of the problem presented in the above model, it is necessary to control the uncertain parameter in the objective function and constraints. Therefore, assuming the parameter as the minimum degree of constraint satisfaction, the controlled model is as follows:

$$Min\ Z = c^t x \quad (19)$$

$s.t.:$

$$a_{ij}x \geq (1 - \alpha')E_1^{q_{ij}} + \alpha' E_2^{q_{ij}}, \quad i \in m, j \quad (20)$$

$$x \geq 0, \quad \alpha \in [0,1] \quad (21)$$

In the above equation, $E_1^{q_{ij}}, E_2^{q_{ij}}$ is the expected value of the fuzzy number of the used demand parameter, which is calculated as follows:

$$E_1^{q_{ij}} = \frac{q_{ij}^1 + q_{ij}^2}{2} \quad (22)$$

$$E_2^{q_{ij}} = \frac{q_{ij}^3 + q_{ij}^4}{2} \quad (23)$$

Based on this, the following relation can be used instead of the $\tilde{q}_{ij}$ parameter:

$$\tilde{q}_{ij} \cong (1 - \alpha') \frac{q_{ij}^1 + q_{ij}^2}{2} + \alpha' \frac{q_{ij}^3 + q_{ij}^4}{2} \quad (24)$$

method is a well-known approach for solving multi-objective optimization problems. By converting all but one of the objective functions into constraints at each stage, the problem is reduced to a single-objective linear programming problem, which can be solved using standard linear programming techniques.

$$Min\ Z_1(x)$$
$$Z_i(x) \leq \varepsilon_i, \quad \forall i = 2,3,\ldots,n \quad (25)$$
$$x \in X$$

In this problem, $X$ represents the feasible solution space, and $Z_n(x)$ is the $n$th objective function of the multi-objective optimization problem. In this method, we reduce the complexity of the objective space and increase the complexity of the solution space, which leads to adding more constraints to the problem and making it more complex compared to the original problem. Another issue with this method is the proper estimation of epsilon $\varepsilon_i$, as an incorrect estimation might result in no feasible solution for the problem. However, with proper adjustment of $\varepsilon_i$, the solution obtained from this method is generally better than that from other methods. Note that the computational complexity of this method is higher than that of the other methods introduced.

The steps of the epsilon-constraint method are as follows:

- One of the objective functions is selected as the main objective function.
- The problem is solved each time with respect to one of the objective functions, and then the optimal values of each objective function are obtained.
- The interval between the two optimal values of the secondary objective functions is divided into a pre-specified number of segments, and a value table for $\varepsilon_2 \ldots \varepsilon_n$ is obtained.
- The problem is solved each time with the main objective function and each of the $\varepsilon_2 \ldots \varepsilon_n$ values.
- The Pareto solutions found are presented.

**NSGA-II Algorithm**

The NSGA-II (Non-dominated Sorting Genetic Algorithm II), initially introduced by Deb et al. (2002), was developed to overcome the limitations of conventional optimization techniques, including computational complexity, absence of elitism, and the necessity of defining a niche parameter. This algorithm employs an elitist strategy to systematically refine the Pareto front, preserving superior solutions from previous generations while utilizing genetic operators to evolve new ones. By integrating selection, crossover, and mutation processes, NSGA-II ensures a well-distributed set of optimal solutions.

In NSGA-II, solutions are first ranked based on their dominance relations and then sorted using a crowding distance metric to maintain diversity (Ahmadianfar et al., 2017). Important tuning parameters such as iteration limits, population size, and mutation and crossover rates are generally determined through empirical testing.

**Execution Steps for NSGA-II:**
1. Generate an initial random population of size.
2. Sort solutions using non-dominated sorting.
3. Assign ranks based on Pareto dominance and apply variation operators to create offspring.
4. Merge parent and offspring populations and categorize solutions using non-dominated ranking.
5. Select individuals for the next iteration using crowding distance.
6. Repeat until convergence criteria are satisfied.

**MOPSO Algorithm**

The Multi-Objective Particle Swarm Optimization (MOPSO) technique, as proposed by Coello et al. (2002), integrates swarm intelligence with an external archive that stores non-dominated solutions approximating the Pareto front. Each particle within the swarm adjusts its trajectory based on individual memory, local best-known solutions, and the best-known global solution discovered by any swarm member.

MOPSO ensures an effective balance between exploration and exploitation by leveraging these guiding principles. Parameter selection, including swarm size, iteration count, mutation rate, and archive size, is typically optimized through experimental tuning.

**MOPSO Execution Process:**

1. Initialize the swarm and identify non-dominated solutions.
2. Store superior solutions in an external archive.
3. Discretize the search space and select reference leaders for swarm guidance.
4. Update memory for each particle.
5. Add newly found non-dominated solutions to the archive and remove dominated entries.
6. Control archive capacity by eliminating redundant solutions.
7. Repeat until the stopping condition is met.

**MOWOA Algorithm**

The Multi-Objective Whale Optimization Algorithm (MOWOA) draws inspiration from the cooperative hunting behavior of humpback whales, particularly their bubble-net feeding technique. This nature-inspired metaheuristic method combines exploration (searching for global solutions) with exploitation (refining discovered solutions) to generate a well-balanced Pareto front.

MOWOA dynamically adjusts search agent positions by leveraging adaptive encircling, spiral movement, and stochastic exploration techniques.

**Key Computational Phases in MOWOA:**

1. Randomly initialize search agents.
2. Assess and record the best-performing solutions.
3. Update positions using a logarithmic spiral movement model.
4. Apply adaptive encircling behaviors to improve convergence.
5. Introduce randomized exploration mechanisms to enhance diversity.
6. Continue iterations until the convergence condition is met.

Table (1) shows the proposed parameter values for each level in metaheuristic algorithms.

**Table 1- Suggested parameter values of meta-heuristic algorithms**

| Optima Level | Level 3 | Level 2 | Level 1 | Parametrs | Algorithm |
|---|---|---|---|---|---|
| Level 3 | 125 | 100 | 75 | $Max\ it$ | |
| Level 2 | 120 | 100 | 50 | $N\ pop$ | NSGA-II |
| Level 1 | 0.07 | 0.06 | 0.05 | $P_c$ | |
| Level 3 | 0.9 | 0.8 | 0.7 | $P_m$ | |
| Level 3 | 125 | 100 | 75 | $Max\ it$ | |
| Level 2 | 120 | 100 | 50 | $N\ particle$ | |
| Level 1 | 2 | 1.5 | 1 | $c_1$ | MOPSO |
| Level 1 | 2 | 1.5 | 1 | $c_2$ | |
| Level 3 | 0.9 | 0.8 | 0.7 | $w$ | |
| Level 3 | 125 | 100 | 75 | $Max\ it$ | |
| Level 2 | 120 | 100 | 50 | $N\ whale$ | MOWOA |
| Level 2 | 3 | 2 | 1 | $A$ | |
| Level 3 | 3 | 2 | 1 | $C$ | |

## Analysis of Results

This section of the paper analyzes various numerical examples for the hub location problem in distribution network issues. Therefore, in this section, the epsilon-constraint method is used for validating and analyzing the sensitivity of the mathematical model, and the NSGA-II, MOPSO, and MOWOA algorithms are used to solve the model in larger sizes.

## Numerical example analysis with epsilon constraint

Initially, to examine the mathematical model, a numerical example with 10 nodes was considered, where the maximum number of hubs created is 3. Due to the lack of real-world data, some of the problem's data are randomly used based on a uniform distribution, as shown in Table (2).

**Table 2- Value parameters according to uniform distribution**

| Value | Parameter | Value | Parameter |
|---|---|---|---|
| 0.6 | $\alpha$ | 250 | $\omega$ |
| 0.8 | $\beta$ | $\sim U(200,500)$ | $f_k$ |
| 1.2 | $EC_{ij}$ | $\sim U(0.1,0.2)$ | $u_k$ |
| 1.3 | $DC_{ij}$ | $\sim U(50,300)$ | $d_{ij}$ |
| 50 | $\varphi$ | $\sim U(2,3)$ | $c_{ij}$ |
| $\sim U(200,300)$ | $\sigma_{ij}$ | $[d_{ij}/10]$ | $\tau_{ij}$ |
| $\sim U(2000,3000)$ | $\vartheta_k$ | 1 | $E_{p1}^{LTO}$ |
| 3 | $E_{p2}^{LTO}$ | 2 | $R_{p1}^{CCD}$ |
| 0.5 | $R_{p2}^{CCD}$ | | |
| $q_{ij}^1 \sim U(60,70) - q_{ij}^1 \sim U(60,70) - q_{ij}^1 \sim U(60,70) - q_{ij}^1 \sim U(60,70)$ | | | $\tilde{q}_{ij}$ |

After designing the small-sized numerical example, the mathematical model was solved using the epsilon-constraint method. As a result of this analysis, 7 efficient solutions were obtained by each method.

Fig. 3 shows the Pareto front and the set of efficient solutions obtained from solving the small-sized numerical example with the epsilon-constraint.

| Solutions | $Z_1$ | $Z_2$ | $Z_3$ |
|---|---|---|---|
| 1 | 2347536.55 | 98050.54 | 531.60 |
| 2 | 2449664.99 | 102275.68 | 440.40 |
| 3 | 2538709.05 | 95627.62 | 549.60 |
| 4 | 2631678.08 | 108230.95 | 331.20 |
| 5 | 2642442.57 | 107688.65 | 331.20 |
| 6 | 2829071.09 | 115277.44 | 222.00 |
| 7 | 3168316.32 | 123287.86 | 112.80 |

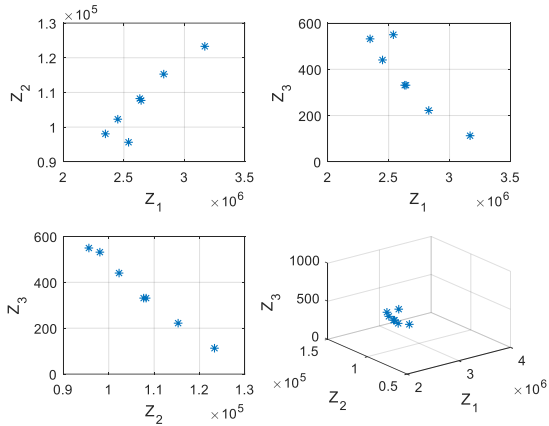

**Figure 3- Pareto front and set of efficient solutions to the small size problem**

The results in Fig. 3 show that 7 efficient solutions were obtained from exact solving methods. Analyzing the efficient solutions reveals that as the penalty costs for exceeding the time window decrease, the total hub location costs increase due to changes in product transfer routes. Also, as the flow transfer time between different nodes decreases, the total hub location network costs also increase.

Examining the first efficient solution to the problem shows that among the three potential nodes, nodes (1) and (5) are selected as main hubs, and other nodes are selected as non-hub nodes. In this analysis, nodes 2-3-4-6-9 are allocated to hub (1), and nodes 7-8-10 are allocated to hub (5). Since direct distribution of goods via freight transport is possible in this model, Table (3) shows all allocations made in the small-sized problem.

Table 3- How to transfer between two nodes in the example of small size

|  | $i_1$ | $i_2$ | $i_3$ | $i_4$ | $i_5$ | $i_6$ | $i_7$ | $i_8$ | $i_9$ | $i_{10}$ |
|---|---|---|---|---|---|---|---|---|---|---|
| $i_1$ | - | Direct | Direct | Direct | Direct | Direct | $k_5$ | $k_5$ | Direct | $k_5$ |
| $i_2$ | Direct | - | $k_1$ | $k_1$ | $k_1$ | Direct | $k_1 \to k_5$ | $k_1 \to k_5$ | $k_1$ | Direct |
| $i_3$ | Direct | $k_1$ | - | Direct | Direct | Direct | $k_1 \to k_5$ | Direct | $k_1$ | Direct |
| $i_4$ | Direct | $k_1$ | Direct | - | $k_1$ | Direct | Direct | Direct | Direct | Direct |
| $i_5$ | Direct | $k_1$ | $k_1$ | $k_1$ | - | $k_1$ | Direct | Direct | $k_1$ | Direct |
| $i_6$ | Direct | Direct | Direct | Direct | $k_1$ | - | $k_1 \to k_5$ | $k_1 \to k_5$ | Direct | Direct |
| $i_7$ | $k_5$ | Direct | $k_5 \to k_1$ | Direct | Direct | $k_5 \to k_1$ | - | $k_5$ | $k_5 \to k_1$ | $k_5$ |
| $i_8$ | $k_5$ | $k_5 \to k_1$ | Direct | Direct | Direct | $k_5 \to k_1$ | $k_5$ | - | $k_5 \to k_1$ | Direct |
| $i_9$ | Direct | $k_1$ | $k_1$ | Direct | $k_1$ | Direct | $k_1 \to k_5$ | $k_1 \to k_5$ | - | Direct |
| $i_{10}$ | $k_5$ | Direct | Direct | Direct | Direct | Direct | $k_5$ | Direct | Direct | - |

In the table above, it is observed that the connection between two nodes, either through one hub, two hubs, or directly, is indicated by the symbol "Direct". Thus, the examined model is valid. Subsequently, the sensitivity of the problem to changes in important parameters is analyzed. Several important parameters were selected and changed, and efficient solution number 1 was compared. Initially, sensitivity analysis was conducted under changes in coefficients α and β. Table (4) shows the variations in the objective functions of the first efficient solution under different coefficients α and β.

Table 4- Changes of the objective functions of the first effective solution in different α and β value

| $Z_3$ | $Z_2$ | $Z_1$ | β | α |
|---|---|---|---|---|
| 532.80 | 97534.66 | 2283599.30 | 0.8 | 0.4 |
| 531.60 | 97826.15 | 2313547.22 | 0.8 | 0.5 |
| 531.60 | 98050.54 | 2347536.55 | 0.8 | 0.6 |
| 531.60 | 98234.21 | 2365899.27 | 0.8 | 0.7 |
| 531.60 | 98234.21 | 2381437.52 | 0.8 | 0.8 |
| 531.60 | 98234.21 | 2056742.16 | 0.6 | 0.6 |
| 531.60 | 98234.21 | 2201533.17 | 0.7 | 0.6 |
| 531.60 | 98050.54 | 2347536.55 | 0.8 | 0.6 |
| 532.80 | 97621.27 | 2467984.81 | 0.9 | 0.6 |
| 534.00 | 97621.27 | 2594981.74 | 1 | 0.6 |

According to the results in Table (4), with increasing coefficients α and β, the network costs increase due to the rise in economic discount factors. Additionally, with an increase in coefficient α, the second objective function value increases, and the third objective function value decreases. In contrast, with an increase in coefficient β, this trend is reversed, with the second objective function value decreasing and the third objective function value increasing. This difference is due to the change in the impact of discount factors resulting from the distribution of transport costs among hubs and non-hub nodes. Fig. 4 shows the variations in the objective function values of the problem under different coefficients α and β.

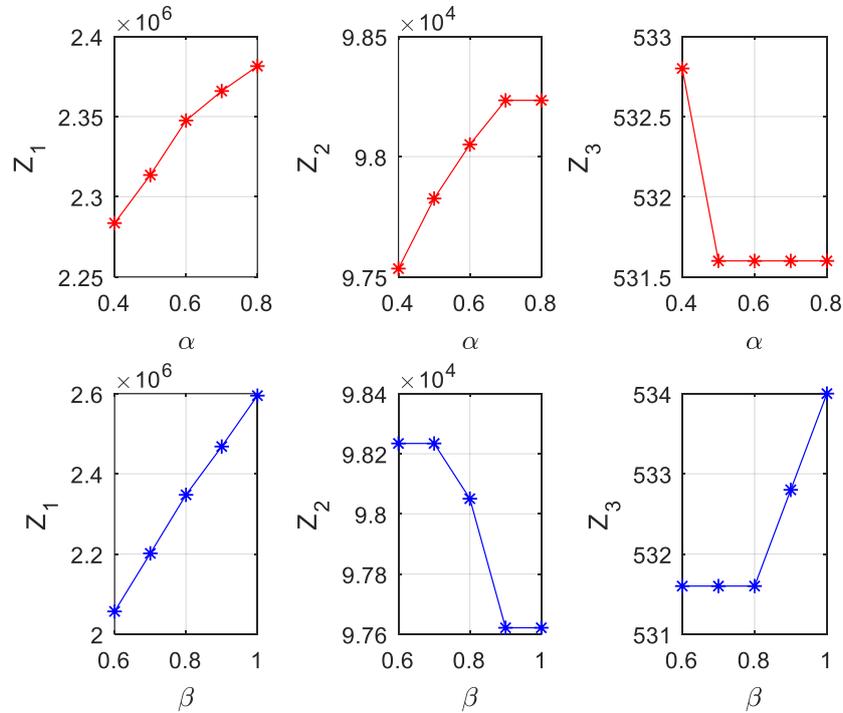

Figure 4- Changes in the values of the objective functions of the problem in different α and β value

Aircraft capacity is another parameter that affects the objective function values of the problem. Table (5) shows the variations in the objective function values for different aircraft capacities ranging from 30 to 70 units.

Table 5- Changes of the objective functions of the first efficient solution in different capacities of the aircraft

| $Z_3$ | $Z_2$ | $Z_1$ | $\varphi$ |
|---|---|---|---|
| 531.60 | 121456.81 | 2347536.55 | 30 |
| 531.60 | 110475.18 | 2347536.55 | 40 |
| 531.60 | 98050.54 | 2347536.55 | 50 |
| 531.60 | 96121.42 | 2347536.55 | 60 |
| 531.60 | 86495.29 | 2347536.55 | 70 |

The analysis of changes in the objective functions at different aircraft capacities for transporting products shows that this parameter only affects the second objective function, and with increasing capacity, the amount of environmentally harmful gas emissions decreases. This is due to the reduction in the number of aircraft used in the hub location network for air transport.

Considering the uncertainty in demand parameters in this paper and the use of fuzzy programming methods to control the mathematical model, changes in the objective function values due to uncertainty rates are shown in Fig. 5. This analysis examines the objective function values of the problem at uncertainty rates of 0.1, 0.3, 0.5, 0.7, and 0.9.

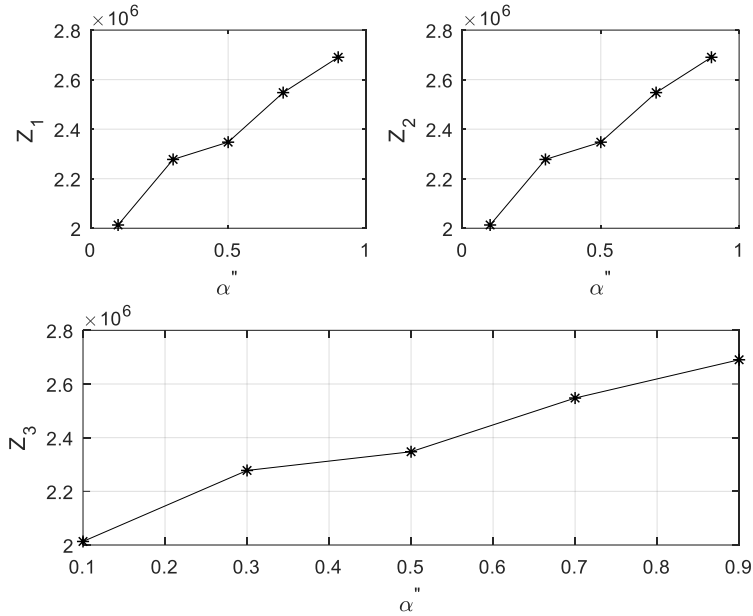

**Figure 5- Changes in the values of the objective functions of the problem in the uncertainty rate**

Fig. 5 shows that with the increase in uncertainty rates, the demand values for nodes have increased. This increase in demand has led to higher product transfer costs through the air transport system and also increased the level of environmental pollutant emissions.

Before analyzing numerical examples of larger sizes, this section addresses solving a small-sized numerical example using NSGA-II, MOPSO, and MOWOA algorithms. The data and size of the numerical example are consistent with the previous example, which was solved using exact methods. Due to the large number of efficient solutions obtained through these algorithms, a comparison of comparison indices and Pareto fronts between different solution methods is provided as shown in Fig. 6.

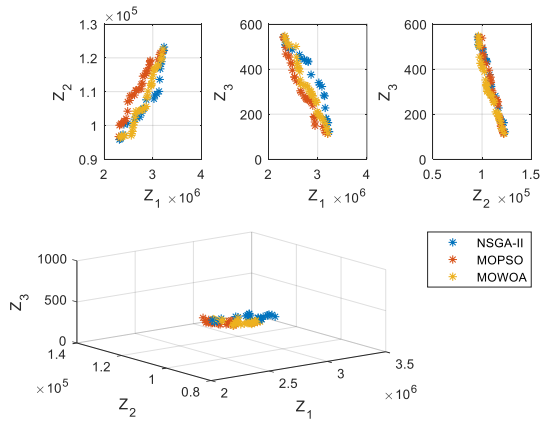

| MOWOA | MOPSO | NSGA II | EC | Factor |
|---|---|---|---|---|
| 38 | 32 | 27 | 7 | NPF |
| 3795.2 | 3594.6 | 3658.6 | 4356.4 | MSI |
| 0.63 | 062 | 0.59 | 0.42 | SM |
| 36.4 | 28.9 | 31.2 | 234.1 | CPT |

**Figure 6- Pareto front comparison indices obtained from solving a small numerical example with meta-heuristic algorithms**

The results of solving the small-sized model with metaheuristic algorithms show that the convergence of the obtained solutions is close to the epsilon-constraint method. Additionally, the solving time using metaheuristic algorithms is significantly less than that of the epsilon-constraint method. In these results, MOWOA achieved the highest NPF and MSI values, NSGA-II achieved the lowest SM value, and MOPSO achieved the lowest CPT value.

## Analysis of numerical examples with meta-heuristic algorithms

After analyzing the small-sized numerical example and the sensitivity analysis of the mathematical model, the analysis of larger-sized numerical examples with different solution methods was conducted. Accordingly, 10 numerical examples were designed as shown in Table (6).

**Table 6- Size of numerical examples in larger size**

| 10 | 9 | 8 | 7 | 6 | 5 | 4 | 3 | 2 | 1 | |
|---|---|---|---|---|---|---|---|---|---|---|
| 150 | 120 | 110 | 100 | 90 | 80 | 60 | 45 | 30 | 15 | # Node |
| 75 | 60 | 50 | 45 | 40 | 35 | 20 | 15 | 10 | 6 | # Hub |

After designing the numerical examples, they were solved using various solution methods. According to the obtained results, only numerical example number 1 was solved by exact methods in less than 1000 seconds. Therefore, Table (7) shows only the comparison indices of efficient solutions in large-sized numerical examples between metaheuristic algorithms.

**Table 7- Indicators of the comparison of efficient solutions in different numerical examples**

| MOWOA | | | | MOPSO | | | | NSGA-II | | | | |
|---|---|---|---|---|---|---|---|---|---|---|---|---|
| CPT | SM | MSI | NPF | CPT | SM | MSI | NPF | CPT | SM | MSI | NPF | Numerical Example |
| 62.81 | 0.24 | 3822 | 50 | 68.19 | 0.11 | 4749 | 38 | 54.60 | 0.42 | 4028 | 32 | 1 |
| 72.74 | 0.16 | 2575 | 48 | 80.01 | 0.23 | 2004 | 19 | 62.03 | 0.48 | 2867 | 36 | 2 |
| 87.20 | 0.33 | 4215 | 41 | 97.17 | 0.27 | 3383 | 40 | 69.70 | 0.23 | 4015 | 48 | 3 |
| 103.0 | 0.20 | 2728 | 35 | 116.3 | 0.21 | 3270 | 37 | 77.74 | 0.37 | 4085 | 43 | 4 |
| 120.1 | 0.12 | 4752 | 40 | 137.2 | 0.18 | 3387 | 15 | 88.43 | 0.28 | 2203 | 32 | 5 |
| 142.1 | 0.40 | 2807 | 43 | 164.4 | 0.43 | 4314 | 37 | 98.49 | 0.43 | 2764 | 40 | 6 |
| 174.7 | 0.20 | 4296 | 37 | 204.7 | 0.27 | 2964 | 21 | 112.9 | 0.41 | 2672 | 34 | 7 |
| 198.7 | 0.28 | 2565 | 43 | 235.6 | 0.46 | 4352 | 20 | 126.7 | 0.17 | 4003 | 40 | 8 |
| 228.5 | 0.38 | 2862 | 41 | 274.2 | 0.26 | 3410 | 29 | 139.3 | 0.44 | 4533 | 35 | 9 |
| 263.3 | 0.24 | 2273 | 38 | 319.7 | 0.41 | 2102 | 32 | 156.5 | 0.50 | 3033 | 34 | 10 |
| 145.3 | 0.255 | 3289 | **41.6** | 169.7 | **0.283** | 3393 | 28.8 | **98.6** | 0.373 | **3420** | 37.4 | Average |

Reviewing the comparison indices indicates that NSGA-II achieved the highest average MSI and the lowest CPT value; MOPSO achieved the lowest average SM, and MOWOA achieved the highest average NPF. Additionally, examining the trend of changes in average CPT shows that with the increase in problem size, CPT increases exponentially. Fig. 7 shows a comparison of average indices in solving large-sized problems with metaheuristic algorithms.

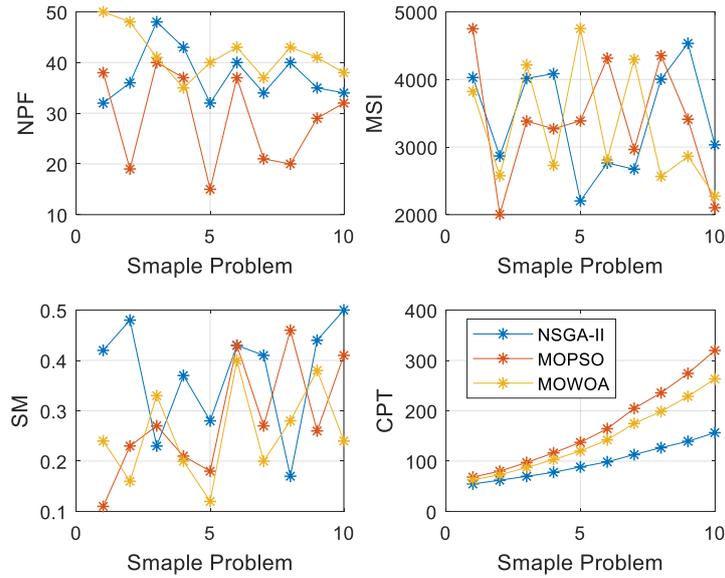

**Figure 7- Comparison of factors in solving large size problems with meta-heuristic algorithms**

The review results indicate that each algorithm was efficient in achieving a specific index, so TOPSIS was used to prioritize the efficient algorithms. In this method, the weight of each index is considered to be 0.25. Table (8) summarizes the ranking results of algorithms using the TOPSIS method.

**Table 8- Ranking meta-heuristic algorithms with TOPSIS**

| Rank | CI | CPT | SM | MSI | NPF | Algorithm |
| --- | --- | --- | --- | --- | --- | --- |
| 1 | 0.770 | **98.6** | 0.373 | **3420** | 37.4 | NSGA-II |
| 2 | 0.666 | 169.7 | **0.283** | 3393 | 28.8 | MOPSO |
| 3 | 0.303 | 145.3 | 0.255 | 3289 | **41.6** | MOWOA |

By comparing the CI index, it can be stated that NSGA-II exhibits higher efficiency in solving the hub location problem for air cargo transport with respect to stability compared to MOPSO and MOWOA.

**Summary and Conclusion**

The significance of this study lies in its comprehensive approach to addressing some of the most pressing challenges in air cargo transportation. With the growing importance of sustainable practices in logistics, the proposed model provides a practical solution for minimizing environmental impacts, which is essential in an industry that contributes significantly to global CO2 emissions. Additionally, by incorporating time-sensitive logistics into the hub network design, the model ensures that service levels are maintained, thereby helping companies to meet customer expectations and avoid penalties for late deliveries.

Moreover, the use of fuzzy logic to account for uncertainty in demand adds a layer of realism to the model, making it more applicable to real-world scenarios where fluctuations in demand are common. This study's findings can help logistics companies optimize their air cargo networks while simultaneously meeting environmental regulations and improving operational efficiency. This paper addressed the modeling and solving of a hub location problem for air cargo transport considering sustainability under uncertainty. The primary goal of establishing and designing a hub network is to exploit cost economies of scale. However, in recent years, the issue of environmental pollutant emissions has become increasingly important, and research in this area is growing. Considering that the air transport industry can have a significant impact on pollutant emissions, this paper aimed to minimize environmental pollutant emissions alongside economic and social aspects. The proposed model seeks to determine nodes as hubs and allocate non-hub nodes to hubs under demand uncertainty.

The results of the mathematical model analysis using fuzzy programming methods showed that as penalty costs for delays decrease, the costs associated with hub location and allocation increase, which also leads to increased environmental pollutant emissions. Specifically, 7 efficient solutions were obtained using the epsilon-constraint method, 27 efficient solutions using NSGA-II, 32 efficient solutions using MOPSO, and 38 efficient solutions using MOWOA. The analyses showed a high convergence of metaheuristic algorithms in achieving efficient solutions with much lower CPT and higher NPF. Analyzing the first efficient solution and examining uncertainty rates revealed that increased demand raises the costs associated with product distribution through the air transport system, and the fixed aircraft capacity relative to the increase in demand leads to greater use of the air transport system, which increases environmental pollutant emissions.

Further analysis with metaheuristic algorithms showed that NSGA-II generally has higher efficiency in achieving CPT and MSI. Meanwhile, MOPSO achieved the best SM, and MOWOA achieved the best NPF. By comparing the CI index, NSGA-II proved to be an efficient algorithm for solving the proposed model in this paper. Future research is recommended to consider different types of aircraft with varying capacities and limited numbers for sending shipments. Additionally, the mathematical model should be considered in a multi-period and dynamic manner. The use of other solving methods to improve efficiency and solve the model in less time and for larger dimensions is also suggested as another recommendation.

**Conflict of Interest**

The authors declare that they have no known competing financial interests or personal relationships that could have appeared to influence the work reported in this paper.